\documentclass[12pt, reqno]{amsart}
\usepackage{graphicx}
\usepackage{amsmath, amsthm, amscd, amsfonts, amssymb, graphicx, color}
\usepackage[bookmarksnumbered, colorlinks, plainpages]{hyperref}
 \newtheorem{theorem}{Theorem}[section]
\newtheorem{lemma}[theorem]{Lemma}

\newtheorem{corollary}[theorem]{Corollary}
\theoremstyle{definition}
\newtheorem{definition}[theorem]{Definition}

\newtheorem{remark}[theorem]{Remark}
\numberwithin{equation}{section}

\begin{document}

%
%
%
%
%
%
%
%
%

\title[K-FRAMES IN HILBERT MODULES]
 {K-FRAMES IN HILBERT $C^{*}$-MODULES}

\author[Abbaspour]{Gh. Abbaspour Tabadkan}

\address{%
Department of mathematics\\
 School of mathematics and computer science\\
Damghan university\\
 Damghan, Iran.\\
}

\email{abbaspour@du.ac.ir}

\author [Arefijamaal]{A.A. Arefijamaal}
\address{
Faculty of Mathematics and Computer Sciences\\
Hakim Sabzevari University\\
Sabzevar, IRAN
}
\email{arefijamal@sttu.ac.ir}
\author [Mahmoudieh]{M. Mahmoudieh}
\address{
Department of mathematics\\
 School of mathematics and computer science\\
Damghan university\\
 Damghan, Iran.\\}
 \email{mahmoudieh@du.ac.ir}

\date{August 27, 2017}

\begin{abstract}
In this paper, firstly we investigate conditions under which the action of an operator on a $K$-frame, remain again a $K$-frame for Hilbert module E. We also give a generalization of Douglas Theorem and we shall use it to prove the sum of two $K$-frame under certain condition is again a $K$-frame. Finally, we characterize the $K$-frame generators for a unitary system in terms of operators.
\end{abstract}

\maketitle


\section{Introduction}
\label{sec:intro}
Frames were first introduced in 1952 by Duffin and Schaeffer \cite{dufin}. Frames can be viewed as redundant bases which are generalization of orthonormal bases.
Many generalizations of frames were introduced, e.g., frames of subspaces \cite{Frames of subspaces}, Pseudo-frames \cite{Pseudo-frame}, G-frames \cite{G-frames}, 
and  fusion frames \cite{Fusion frames}. Recently, L. Gavruta introduced the concept of {$ K\text{-frame} $} for a given bounded operator K on Hilbert space in \cite{Gavruta}. Hilbert $C^{*}$-modules arose as generalizations of the notion ’Hilbert space’. The basic idea was to consider modules over $C^{*}$-algebras instead of linear spaces and to allow the inner product to take values in
the $C^{*}$-algebra of coefficients being $C^{*}$-(anti-)linear in its arguments \cite{lance}. In \cite{frank1} authors generalized frame concept for operators in Hilbert $C^{*}$-modules. The paper is organized as follows.
In Section 2, some notations and preliminary results of Hilbert Modules, their frames and $K\text{-frames}$ are given. In Section 3, we study the action of operators on $K$-frames and under certain conditions, we shall show that it is again a $K$-frame.
The next section is devoted to sum of $K$-frames, to show that the sum of two $K$-frames under certain conditions is again a $K$-frame we need to say a generalization of the Douglas Theorem \cite{Jing_Yu}, which may interest by its own. Finally, in the last section, we consider a unitary system of operators and characterize the $K$-frame generators in terms of operators. We also look forward to sum of two $K$-frame generators to be a $K$-frame generator.

\section{Preliminaries}
In this section we give preliminaries about frames, $ K $-frames for Hilbert space and Hilbert module and related operators which we need in the sections following.
A finite or countable sequence $ \{f_{k}\} _{{k \in \mathbb{J}}} $ is called a frame for a separable Hilbert space $ H $ if there exist constants $ A , B > 0 $ such that
\[A{\left\| {{\rm{ }}f} \right\|^2} \le \sum\limits_{k = 1}^\infty  {{{\left| {\left\langle {f,{f_k}} \right\rangle } \right|}^2}}  \le B{\left\| {{\rm{ }}f} \right\|^2},\quad {\rm{      }}\forall {\rm{f}} \in H.\]

Frank and Larson \cite {frank1} introduced the notion of frames in Hilbert {$ C^{*}\text{-modules} $} as a generalization of frames in
Hilbert spaces. 
A (left) Hilbert $C^{*}$-module over the $C^{*}$-algebra $\mathcal{A}$
is a left $\mathcal{A}$-module $E$ equipped with an $\mathcal{A}$-valued inner product satisfy the following
conditions:
\begin{enumerate}
\item
$\langle x,x\rangle\geq 0$ for every $x\in E$ and $\langle x,x\rangle =0$ if and only if $x=0$,
\item
$\langle x,y\rangle =\langle y,x\rangle^*$ for every $x,y\in E$,
\item
$\langle\cdot,\cdot\rangle$ is $\mathcal{A}$-linear in the first argument,
\item
$E$ is complete with respect to the norm $\| x\|^2=\|\langle x,x\rangle\|_\mathcal{A}$.
\end{enumerate}
Given Hilbert $C^{*}$-modules $E$ and $F$, we denote by $ L_{\mathcal{A}}(E,F) $  or $ L(E, F)$ the set of all adjointable operators from $E$ to $F$
 i.e. the set of all maps $T: E\rightarrow F$ such that there exists $T^*: F\rightarrow E $ with the property
$$\langle Tx, y\rangle = \langle x, T^*y\rangle$$ for all $x\in E$, $ y\in F.$
It is well-known that each adjointable operator is necessarily bounded $\mathcal{A}$-linear in the sense $T(ax) = aT(x)$, for all $a\in \mathcal{A}, x\in E$. We denote $ L(E) $ for $ L(E,E) .$ In fact  $ L(E) $ is a \linebreak $C^{*}$-algebra.
%
Let $\mathcal{A}$ be a $C^{*}$-algebra and consider
$$\ell^{2}(\mathcal{A}) := \{\{a_n\}_n\subseteq \mathcal{A}:~\sum_{n}a_n a^{*}_{n}\text{\qquad converges in norm in}~\mathcal{A} \}.$$
It is easy too see that $\ell^2(\mathcal{A})$ with pointwise operations and the inner product
$$\langle\{a_n\}, \{b_n\}\rangle=\sum_{n}a_n b^{*}_{n},$$
becomes a Hilbert $C^{*}$-module which is called the standard Hilbert {$ C^{*}\text{-module} $} over $\mathcal{A}$.
Throughout this paper, we suppose $ E $ is a Hilbert $\mathcal{A}$-module, $ \mathbb{J} $ a countable index set. Also we denote the range of $T \in L(E) $ by $R(T) $, and kernel of $ T $ by $ N(T) $.
A Hilbert $\mathcal{A}$-module $E$ is called finitely generated (countably generated) if there exist a finite subset $\{x_1, ..., x_n\}$ (countable set $\{x_n\}_{n \in \mathbb{J} }$)
of $E$ such that $E$ equals the closed $\mathcal{A}$-linear hull of this set. The basic theory of
Hilbert $C^{*}$-modules can be found in \cite{lance}.\\ 
The following lemma found the relation between the range of an operator and kernel of its adjoint operator.\\ 
\begin{lemma} \label{wegge} (\cite{wegge-Olsen}, Lemma 15.3.5; \cite{lance}, Theorem 3.2 )  
Let ${ T\in L(E,F)} $, then 
\begin{enumerate}
\item
$ N(T)=N(\vert T \vert) $, $ N({T}^*) ={R(T)}^\perp  $,  $ {N({T}^*)}^\perp ={R(T)}^{\perp\perp} \supseteq \overline{R(T)}$;\
\item
$ R(T) $ \text{ is closed if and only if} $ R({T}^*) $ \text{ is closed, and in this case}  $ R(T) $ and $ R({T}^*) $\text{are orthogonally complemented with}  $ R(T)=$ \\ ${N({T}^*)^\perp} $ \text{and} $ R({T}^*)= {N(T)}^\perp .$
\end{enumerate}
\end{lemma}

The following theorem is extended Douglas theorem \cite{Doglas} for Hilbert modules.
\begin{theorem} \label{dog_md} \cite{Jing_Yu}
Let  $ {T}^{\prime }\in L(G,F) $ and $ T\in L(E,F) $ with $ \overline{R(T^{*})} $ orthogonally complemented. The following statements are equivalent: 
\begin{enumerate}
\item
 $   {T}^{\prime} {T}^{\prime *}\leq \lambda T{T}^* $ for some $ \lambda  > 0  $;
\item
There exists $ \mu > 0 $ such that $ \Vert {T}^{'*}z \Vert \leq \mu \Vert {T}^{* }z \Vert $ for all $ z \in F;$
\item
There exists $ D \in L(G,E)  $such that $  {T}^{\prime }=TD $, i.e. the equation $ TX= {T}^{\prime } $ has a solution;
\item $ R( {T}^{\prime }) \subseteq R(T) $.
\end{enumerate}

\end{theorem}

Here, we recall the concept of frame in Hilbert $C^{*}$-modules which is defined in \cite{frank1}.
Let $ E $ be a countably generated Hilbert module over a unital $C^{*}$-algebra $\mathcal{A}$. A sequence $\{x_n\}\subset E$ is said
to be a \emph{frame} if there exist two constant $C, D>0$ such that
\begin{equation}\label{frame intro}
C\langle x, x\rangle \leq \sum_{n}\langle x, x_n\rangle\langle x_n, x\rangle\leq D\langle x, x\rangle {\text { for all }} x\in E.
\end{equation}
The optimal constants (i.e. maximal for $C$ and minimal for $D$) are called frame bounds. If the sum in \eqref{frame intro} converges in norm, the frame is called \emph{standard frame}. In this paper all frames consider standard frames.
The sequence $\{x_n\}$ is called a \emph{Bessel sequence} with bound $D$ if the upper inequality in \eqref{frame intro} holds for every $x\in E$. Let $ \{x_{j}\} _{{j \in \mathbb{J} }} $ be a Bessel sequence  for Hilbert module $E$ over $\mathcal{A}$. The operator ${T: E \rightarrow \ell^2(\mathcal{A})} $ defined by $ Tx =  \{ \langle x , x_{j} \rangle \} _{{j \in \mathbb{J}}} $ is called the analysis operator. The adjoint operator  $T^{*}:  \ell^2(\mathcal{A}) \rightarrow E  $ is given by  $$ T^{*} \{c_{j}\} _{{j \in \mathbb{J} }} = \sum_{j \in \mathbb{J}}c_j x_{j},  $$
where is called the \emph{pre-frame operator} or the \emph{synthesis operator}.
By composing $T$ and $T^{*}$, we obtain the \emph{frame operator} $ S: E \rightarrow E $ given by 
$$ Sx = T^{*}Tx = \sum_{j \in \mathbb{J}} \langle x, x_j\rangle x_{j}, \: x \in E. $$
In the case where $ \{x_{j}\} _{{j \in \mathbb{J} }} $ is a frame, the frame operator is positive and invertible, also it is the unique operator in  $ L(E) $ such that the reconstruction formula 
$$ x = \sum_{j \in \mathbb{J}} \langle x, S^{-1}x_j \rangle x_{j} = \sum_{j \in \mathbb{J}} \langle x, x_j \rangle S^{-1} x_{j} , \: x \in E,$$ holds. It is easy to see that the sequence $ \{S^{-1} x_{j}\} _{{j \in \mathbb{J} }} $
is a frame for $E.$ The frame $ \{S^{-1} x_{j}\} _{{j \in \mathbb{J} }} $ is said to be the \emph{canonical dual}  frame of $ \{x_{j}\} _{{j \in \mathbb{J} }}. $
	
	\begin{theorem} \label{NBessel} [ see \cite{Najati}, proposition 2.2 ]
	Let $\{x_{n}\}_{n \in \mathbb{J}} $ be a sequence in $E$ such that $ \sum_{n \in \mathbb{J}} {c_{n}}{x_{n}} $   converges for all  $ c = \{c_{n}\}_{n \in \mathbb{J}} \in {\ell}^{2}(\mathcal{A}) $. Then  $ \{x_{n}\}_{n \in \mathbb{J}} $ is a Bessel sequence in $E.$
   \end{theorem}

	\begin{theorem} \cite{Jing}
	Let $E$ be a finitely or countably generated Hilbert module
	 over a unital $C^{*}$-algebra $\mathcal{A}$, and  $ \{x_{n}\} _{{n \in \mathbb{J}}} $ be a sequence in $E$. Then  $ \{x_{n}\} _{{n \in \mathbb{J}}} $ is a
	frame for $E$ with bounds $ C$ and $D $ if and only if
	$$ C {\Vert x \Vert}^{2} \leq \Vert \sum_{n}\langle x, x_n\rangle\langle x_n, x\rangle \Vert \leq D {\Vert x \Vert}^{2}. $$
	\end{theorem}
In compare with to $K$-frames on Hilbert space, Najati in \cite{Najati} define atomic system and a $K$-frame on Hilbert module.
	\begin{definition} \label{atomic_def}
	 A sequence $ \{x_{n}\}_{n \in \mathbb{J}} $ of $E$ is called an atomic system for $ K \in L(E) $ if the following statement hold:
	\begin{enumerate}
	\item
	The series $ \sum_{n \in \mathbb{J}} {c_{n}}{x_{n}} $   converges for all  $ c = \{c_{n}\}_{n \in \mathbb{J}} \in {\ell}^{2}(\mathcal{A}) $;
	\item
	There exists $ C > 0 $  such that for every $ x \in E $ there exists  $ \{a_{n,x}\}_{_{n \in \mathbb{J}}} \in {\ell}^{2}(\mathcal{A})$ such that $ \sum_{n \in \mathbb{J}}{a}_{n,x}{{a}^{*}}_{n,x} \leq C\langle x,x\rangle $ and \linebreak $ {Kx = \sum_{n \in \mathbb{J}} {a_{n,x}}{x_{n}}} $. 
	\end{enumerate}
	\end{definition}

By Theorem \ref{NBessel}, the condition (1), in the above definition, actually says that $ \{x_{n}\}_{n \in \mathbb{J}} $ is a Bessel sequence.

	\begin{theorem} \cite{Najati}
	If $K \in L(E)$, then there exists an atomic system for $K$.
	\end{theorem}

\begin{theorem} \cite{Najati}
Let $ \{x_{n}\}_{n \in \mathbb{J}} $ be a Bessel sequence for $E$ and $ {K \in L(E)} $. Suppose that $ T \in L(E,{\ell}^{2}(\mathcal{A})) $ is given by $ T(x) = \{\langle x , x_{n} \rangle \}_{n \in \mathbb{J}} $ and $\overline{R(T)}$ is orthogonally complemented. Then the following statements are equivalent: 
\begin{enumerate}
\item 
$ \{x_{n}\}_{n \in \mathbb{J}} $ is an atomic system for $K$;
\item
There exist $C,B > 0 $ such that 
$$ B {\Vert K^{*}x \Vert}^{2} \leq \Vert \sum_{n}\langle x, x_n\rangle\langle x_n, x\rangle \Vert \leq C {\Vert x \Vert}^{2} ;$$
\item
There exist $ D \in L(E,{\ell}^{2}(\mathcal{A})) $ such that $ K = T^{*}D.$
\end{enumerate}
\end{theorem}
	\begin{definition}
	Let $E$ be a Hilbert $\mathcal{A}$-module, $ \{x_{n}\}_{n \in \mathbb{J}} \subset E$ and ${ K \in L(E) } $. The sequence $ \{x_{n}\}_{n \in \mathbb{J}} $ is said to be a $K$-frame if there exist constant $ C, D > 0 $ such that 
	\begin{equation}\label{kframe def}
	C\langle K^{*}x, K^{*}x\rangle \leq \sum_{n}\langle x, x_n\rangle\langle x_n, x\rangle\leq D\langle x, x\rangle , \: \: x \in E.
	\end{equation}
	\end{definition}
The following theorem gives a characterization of $K$-frames using linear adjiontable operators.
 \begin{theorem} \label{L oper} \cite{Najati}
  Let ${ K \in L(E)}$ and $ \{x_{n}\}_{n \in \mathbb{J}} $ be a Bessel sequence for $E$. Suppose that $ T \in L( E , {\ell}^{2}(\mathcal{A}))$ is given by $ T(x) = \lbrace \langle x , x_{n} \rangle \rbrace _{n \in \mathbb{J}} $ and $\overline{R(T)}$ is orthogonally complemented. Then $ \{x_{n}\}_{n \in \mathbb{J}} $ is a $K$-frame for $E$ if and only if there exist a linear bounded operator $ L :  {\ell}^{2}(\mathcal{A})  \rightarrow E $ such that $ Le_{n} = x_{n} $ and $ R(K) \subseteq R(L) $, where $ \{e_{n}\}_{n } $ is the canonical orthonormal basis for $ {\ell}^{2}(\mathcal{A})$.
  
 \end{theorem}
\section{Operators On $K$-frames } 
In this section we study the action of an operator on a $ K $-frame. The following lemma shows that the action of an adjointable operator on a Bessel sequence is again a Bessel sequence.
\begin{lemma} \label{bessel}
Let $E$ be a Hilbert $\mathcal{A}$-module and $ \{ x_{n} \} _{{n \in \mathbb{J}}} $ be a Bessel sequence, then  $ \{Mx_{n}\} _{{n \in \mathbb{J}}} $ is a Bessel sequence for every $ M \in L({E})$. 
\begin{proof}
Since $ \{ x_{n} \} _{{n \in \mathbb{J}}} $ is a Bessel sequence then there exists constant D such that 
\begin{equation*}
\sum_{n}\langle x, x_n\rangle\langle x_n, x\rangle\leq D\langle x, x\rangle
\end{equation*}
for every $x\in E$.
So 
\begin{equation*}
\begin{split}
\sum_{n}\langle x, Mx_n\rangle  \langle Mx_n, x\rangle
& =\sum_{n}\langle M^{*}x, x_n\rangle\langle x_n, M^{*}x\rangle \\
& \leq D\langle M^{*}x, M^{*}x\rangle\\
& = D\langle MM^{*}x, x\rangle \\
& \leq D \Vert M \Vert ^{2}\langle x, x\rangle \\ 
\end{split}
\end{equation*}
for every $x\in E$. 
\end{proof}
\end{lemma}
\begin{theorem}
Let $E$ be a Hilbert $\mathcal{A}$-module, $ K \in L({E}) $  and $ \{x_{n}\} _{{n \in \mathbb{J}}} $ be a $K$-frame for E. Let $ M \in L({E}) $ with $ R(M) \subset R(K)$  and $\overline{R(K^{*})}$ 
orthogonally complemented. Then $ \{x_{n}\} _{{n \in \mathbb{J}}} $ is an M-frame for E.

\begin{proof}
Since  $ \{x_{n}\} _{{n \in \mathbb{J}}} $ is a $K$-frame then there exist positive numbers $\mu $ and $ \lambda $  such that \\
\begin{align} \label{k-frame}
 \lambda \langle K^{*}x , K^{*}x \rangle \leq \sum_{n}\langle x, x_n\rangle\langle x_n, x\rangle\leq \mu \langle x, x\rangle 
\end{align}
Using the theorem \ref{dog_md} by the fact that $ R(M) \subset R(K)$ shows that, $ MM^{*}\leq \lambda^{'} KK^{*} $ for some $\lambda^{'} > 0$. So 
$$ \langle MM^{*}x , x \rangle \leq \lambda^{'} \langle KK^{*}x , x \rangle $$
So
$$  \dfrac{{\lambda}}{{\lambda}^{'}}  \langle MM^{*}x , x \rangle \leq \lambda \langle K^{*}x , K^{*}x \rangle $$ 
From \ref{k-frame}, we have
$$ {\dfrac{{\lambda}}{{\lambda}^{'}} }  \langle MM^{*}x , x \rangle \leq \sum_{n}\langle x, x_n\rangle\langle x_n, x\rangle\leq \mu \langle x, x\rangle. $$ \\
Therefor $ \{x_{n}\} _{{n \in \mathbb{J}}} $ is an $M$-frame with bound $ \dfrac{{\lambda}}{{\lambda}^{'}} $ and $ \mu $ for $ E.$
\end{proof}
\end{theorem}

In the following theorem we obtain the result of the last theorem by different conditions.

\begin{theorem}
Let $ \{x_{n}\} _{{n \in \mathbb{J}}} $ is a $K$-frame for Hilbert $\mathcal{A}$-module E, suppose  $ T \in L(E,l^{2}(A))$ with $ T(x)= \lbrace \langle x , x_{n} \rangle \rbrace $ for every $ x \in E $ ,  $\overline{R(T)} $ orthogonally complemented and $ M \in L(E)$ such that $ R(M) \subset R(K) $. Then $ \{x_{n}\} _{{n \in \mathbb{J}}} $ is an $ M$-frame for $E$.
\begin{proof}
By Theorem \ref{L oper}, there exist $ L : {\ell}^{2}(\mathcal{A}) \rightarrow E $ such that $ {Le_{n} = f_{n}}$, $ R(K) \subset R(L)$. Then $ R(M) \subset R(L)$, now again by Theorem \ref{L oper}, 
 we have $ \{x_{n}\} _{{n \in \mathbb{J}}} $ is an $M$-frame for $E$.
\end{proof}
\end{theorem}

\begin{theorem} \label{surjective}
Let $ E $ be a Hilbert $\mathcal{A}$-module and $ K \in L(E)$ with the dense range. Let $ \{x_{n}\} _{{n \in \mathbb{J}}} $ be a $K$-frame for E and $ T \in L(E)$ has closed range. If $ \{Tx_{n}\} _{{n \in \mathbb{J}}} $ is a $K$-frame for $ E $, then $ T $ is surjective.
\begin{proof}
Suppose that $K^{*}x = 0  $ for $ x \in E$, then for each $ y \in E $, $ { \langle Ky, x \rangle = \langle y, K^{*}x \rangle = 0} $. So $ \langle z, x \rangle = 0 $ for each $ z \in E $, since $ R(K) $ is dense in $ E $. Thus $ x = 0 $ and $ K^{*} $ is injective.	
We shall show that $ T^{*}$ is injective too. Note that 
$ \{Tx_{n}\} _{{n \in \mathbb{J}}} $ is a $K$-frame for E with  bounds $ \lambda $ and $  \mu ,$  hence 
$$ \lambda \Vert K^{*}x \Vert^{2} \leq \Vert \sum_{n}\langle x, Tx_n\rangle\langle Tx_n, x\rangle \Vert \leq \mu \Vert x \Vert^{2}.$$
for $ T^{*}x \in E $ and therefore,
$$ \lambda \Vert K^{*}x \Vert^{2} \leq \Vert \sum_{n}\langle T^{*}x, x_n\rangle\langle x_n, T^{*}x\rangle \Vert \leq \mu \Vert x \Vert^{2}.$$
If $ x \in N(T^{*}) $ then $ T^{*}x = 0 $ so $ \langle T^{*}x , x_{n} \rangle = 0 $ for each $ n \in \mathbb{N}  $, and 
 so $ K^{*}x = 0 $ by the last inequality. On the other hand 
 $ K^{*} $ is injective, so $x = 0 $, and so $ T^{*} $ is injective. Therefore $$  E = N(T^{*}) + \overline{R(T)} = \overline{R(T)} = R(T),  $$ and this complete the proof.
\end{proof}
\end{theorem}

\begin{theorem}{\label{rang_T}}
	Let $ K \in L(E) $ and $ \{x_{n}\} _{{n \in \mathbb{J}}} $ be a $K$-frame for E. If  $ T \in L(E) $ with closed range such that 
	 $ \overline{ R(TK)} $ is orthogonal complemented and $ KT = TK $. Then $ \{Tx_{n}\} _{{n \in \mathbb{J}}} $ is  a $K$-frame for $ R(T). $
\begin{proof}
 Xu and Sheng in \cite{sheng} show that if $ T $ has closed range then $ T $ has Moore-Penrose inverse operator $ T^{\dagger} $ such that $ TT^{\dagger}T = T $ and $ T^{\dagger}TT^{\dagger} = T^{\dagger} $. So $ TT^{\dagger} \vert_{R(T)} = I_{R(T)} $ and ${ (TT^{\dagger})^{*} = I^{*} = I =  TT^{\dagger}.} $ For every $ x \in R(T) $ we have 
\begin{equation*}
 \begin{split}
    \langle K^{*}x , K^{*}x \rangle & = \langle (TT^{\dagger})^{*} K^{*}x , (TT^{\dagger})^{*} K^{*}x \rangle \\
    & = \langle {T^{\dagger}}^{*}T^{*} K^{*}x , {T^{\dagger}}^{*}T^{*} K^{*}x \rangle \\
    & \leq \Vert {(T^{\dagger})}^{*} \Vert^{2} \langle T^{*} K^{*}x , T^{*} K^{*}x \rangle 
 \end{split}
\end{equation*}
and so\\ 
$$ \Vert (T^{\dagger})^{*} \Vert^{-2} \langle K^{*}x , K^{*}x \rangle \leq \langle T^{*} K^{*}x , T^{*} K^{*}x \rangle. $$ \\
Since $ \{x_{n}\} _{{n \in \mathbb{J}}} $ is a $K$-frame and $ R(T^{*}K^{*}) \subset R(K^{*}T^{*}) $, if $ \lambda $ is a lower bound then by using Theorem \ref{dog_md}, there exists some $ \lambda^{'} > 0 $ such that \\
\begin{equation*}
  \begin{split}
    \sum_{n}\langle x, Tx_n\rangle\langle Tx_n, x\rangle & = \sum_{n}\langle T^{*}x, x_n\rangle\langle x_n, T^{*}x\rangle \\
    & \geqslant \lambda \langle  K^{*}T^{*}x ,  K^{*}T^{*}x \rangle \\
    & \geqslant \lambda^{'} \lambda \langle T^{*} K^{*}x , T^{*} K^{*}x \rangle \\
    & \geqslant \lambda^{'} \lambda \Vert (T^{\dagger})^{*} \Vert^{2} \langle K^{*}x , K^{*}x \rangle.
  \end{split}
\end{equation*}
This is the lower inequality for $ \{Tx_{n}\} _{{n \in \mathbb{J}}} $. On the other hand by Lemma \ref{bessel},  $ \{Tx_{n}\} _{{n \in \mathbb{J}}} $ is  a Bessel sequence, so $ \{Tx_{n}\} _{{n \in \mathbb{J}}} $ is a $K$-frame for Hilbert module $ R(T). $
\end{proof}
\end{theorem}

\begin{theorem}{\label{co_iso}}
	Let E be a Hilbert $\mathcal{A}$-module, $ K \in L(E) $ and $ \{x_{n}\} _{{n \in \mathbb{J}}} $ be a $K$-frame for $ E $, and $ T \in L(E) $ is a co-isometry such that \linebreak ${ R(T^{*}K^{*}) \subset R(K^{*}T^{*})} $ with $ \overline{ R(TK)} $ orthogonal complemented. Then \linebreak $ \{Tx_{n}\} _{{n \in \mathbb{J}}} $ is a $K$-frame for E.
	\begin{proof}
Using Lemma \ref{bessel}  $ \{Tx_{n}\} _{{n \in \mathbb{J}}} $ is a Bessel sequence. By Theorem \ref{dog_md}, there exist $\lambda^{'} > 0 $ such that $ \Vert T^{*}K^{*}x \Vert^{2}  \leq \lambda^{'} \Vert K^{*}T^{*}x \Vert^{2} $, for each $ x \in E $.
	Suppose  $\lambda $ is a lower bound for the $K$-frame  $ \{x_{n}\} _{{n \in \mathbb{J}}} $. Since $ T $ is a co-isometry, then
		
		\begin{equation*}
		\begin{split}
	\dfrac{\lambda}{\lambda^{'}} \Vert K^{*}x \Vert^{2} & = \dfrac{\lambda}{\lambda^{'}} \Vert T^{*}K^{*}x \Vert^{2}\\
		& \leq \lambda \Vert K^{*}T^{*}x \Vert^{2} \\
		& \leq  \sum_{n}\langle T^{*}x, x_n\rangle\langle x_n, T^{*}x\rangle \\
		& = \sum_{n}\langle x, Tx_n\rangle\langle Tx_n, x\rangle .
		\end{split}
		\end{equation*}  
		which implies that $ \{Tx_{n}\} _{{n \in \mathbb{J}}} $ is a $ K $-frame for $ E. $
	\end{proof}
\end{theorem}


\begin{remark}
If $ K \in L(E) $ with dense range, $ T \in L(E) $ with closed range such that $ TK = KT $ and $ \{x_{n}\} _{{n \in \mathbb{J}}} $ is a $K$-frame for $E$. Then  $ \{Tx_{n}\} _{{n \in \mathbb{J}}} $ is a $K$-frame for $E$ if and only if  $T$ is surjective.
\end{remark}

 \begin{theorem}
 	 Let $ K \in L(E) $ with dense range and $ \{x_{n}\} _{{n \in \mathbb{J} }} $ is a $K$-frame for $E$. Let  $ T \in L(E) $ with closed range. If $ \{Tx_{n}\} _{{n \in \mathbb{J} }}$ and $\{T^{*}x_{n}\} _{{n \in \mathbb{J}}} $ are $K$-frames for $E$ then $T$ is invertible.
 \begin{proof}
 	 By Theorem \ref{surjective}, $ T $ is surjective. Since $\{T^{*}x_{n}\} _{{n \in \mathbb{J}}} $ is a $K$-frame for $ E $ then there exist positive numbers  $ \mu $ and $ \lambda $ such that for every $ x \in E $
	 $$ \lambda \Vert K^{*}x \Vert^{2} \leq  \Vert {\sum_{n}\langle x, T^{*}x_n\rangle\langle T^{*}x_n, x\rangle } \Vert \leq \mu \Vert x \Vert^{2} $$ \\
	 So for $ x \in N(T) $ we have 
	 $$ \lambda \Vert K^{*}x \Vert^{2} \leq  \sum_{n}\langle x, T^{*}x_n\rangle\langle T^{*}x_n, x\rangle = 0 $$
	 Then $ \Vert K^{*}x \Vert^{2} = 0 $, so $ x \in N(K^{*}) .$ 
	 On the other hand $ K \in L(E) $ has dense range so $ K^{*} $ is injective and so $ T $ is injective.
 \end{proof} 
 \end{theorem}

%
\section{Sums of $K$-frames}
In this section we shall show that the sum of two $K$-frames in a Hilbert $ C^{*} $-module under certain conditions is again a $K$-frame. It is proved, in Hilbert space case, by Ramu and Johnson \cite{Johnson}.	In the proof of Theorem 3.2 of \cite{lance}  indicates that if $ T $ has closed range then $ R(T^{*}T) $  is closed and $ R(T ) = R(T^{*}T)$. The following theorem says that this result still holds for 
	adjointable operators between Hilbert $ C^{*}$-modules (even though $ \overline{ R(T^{*}) } $ 
	may not be complemented). 

	\begin{theorem} \cite{lance}
	For $ T $ in $ L(E,F) $, the sub-spaces  $ R(T^{*}) $  and $ R(T^{*}T) $ have the same closure. 
	
	\end{theorem}
   In \cite{Sharifi}, Sharifi show that the conversely of the above theorem is also true.
	\begin{theorem}[Lemma 1.1, \cite{Sharifi}] \label{Sharifi}
	Suppose $T \in  L(E ) $, then the operator $T$
	has closed range if and only if  $ R(TT^{*}) $ has closed rang. In this case, $R(T ) = R(TT^{*} )$.
	\end{theorem} 
 
\begin{corollary} \label{tarcsay} 
	Suppose $T \in L(E)^{+}$, 
	then $R(T)$ is closed if and only if $R(T^{1/2})$ is closed. In this case, $R(T) = R(T^{1/2}) $.  
    \begin{proof}
    	The proof is immediately consequence of replacement $T$ by $ T^{1/2} $ in the above theorem.
    \end{proof} 
\end{corollary}

\begin{theorem} \label{sqrt theorem}
Let E be a Hilbert module and $ A, B \in L(E)$ such that $ R(A)+R(B) $ is closed. Then 
$$ R(A)+R(B) = R(({AA^{*} + BB^{*}})^{\frac{1}{2}} ) $$.
\begin{proof}
Define $ T \in L( E \oplus E ) $ by $ T:= \begin{bmatrix}
 A &  B \\ 
0 & 0
\end{bmatrix} $ 
then 
$ T^{*}= 
	\begin{bmatrix}
		 A^{*} & 0 \\ 
		 B^{*} & 0 
	\end{bmatrix} $ and \\
$$ TT^{*}=
	 \begin{bmatrix}
		 A &  B \\ 
		0 & 0 \end{bmatrix} 
 \begin{bmatrix}
 A^{*} & 0 \\ 
 B^{*} & 0 
\end{bmatrix}
=\begin{bmatrix}
 AA^{*}+BB^{•} & 0\\ 
0 & 0 \end{bmatrix}. $$ 

So we have 

$$ (TT^{*})^{1/2}= \begin{bmatrix}
 (AA^{*}+BB^{*})^{1/2} & 0\\ 
0 & 0  \end{bmatrix}.$$  
On the other hand
 $$ T  \begin{bmatrix}
  E  \\
  E  
 \end{bmatrix}
 = \begin{bmatrix}
  A &  B \\ 
0 & 0
\end{bmatrix}  
\begin{bmatrix}
 E  \\
 E  
\end{bmatrix}$$
 thus 
$$ R(T) = R(A) + R(B) \oplus \{ 0 \}. $$
Since $ R(T) = ( R(A) + R(B) ) $ is closed then by Theorem \ref{Sharifi}, $ { R(T) = R(TT^{*}) } $, but by the Corollary \ref{tarcsay},  $ R(TT^{*}) = R((TT^{*})^{1/2})$. So we have 
$$ R(A) + R(B) = R((AA^{*} + BB^{*})^{1/2}). $$
\end{proof} 
\end{theorem}
The following theorem is a generalization of Douglas theorem [Theorem 1.1, \cite{Jing_Yu} ], for Hilbert modules. 

\begin{theorem} \label{dog for sum}

Let $ A,B_{1},B_{2} \in L(E) $ and $ R(B_{1}) + R(B_{2}) $ 
 is 
 closed. The following statements are equivalent.
\begin{enumerate}
\item
$ R(A) \subset R(B_{1}) + R(B_{2}) $;
\item
$ AA^{*} \leq \lambda ( B_{1}{B_{1}}^{*} + B_{2}{B_{2}}^{*} ) $ for some $ \lambda > 0 ;$ 
\item
There exist $ X, Y \in L(E) $ such that $ A = B_{1}X + B_{2}Y $. 
\end{enumerate}
\begin{proof}
$(1) \Longrightarrow (2)$: Suppose $ R(A) \subset R(B_{1}) + R(B_{2}) $ then by Theorem \ref{sqrt theorem}, we have 
\begin{equation*}
\begin{split}
R(A) & \subset R(B_{1}) + R(B_{2})\\
& = R((B_{1}{B_{1}}^{*} + B_{2}{B_{2}}^{*})^{1/2})\\
\end{split}
\end{equation*}
thus by Theorem \ref{dog_md}, $ AA^{*} \leq \lambda ( B_{1}{B_{1}}^{*} + B_{2}{B_{2}}^{*})$ for some $ \lambda > 0 $.
 \\
 $(2) \Longrightarrow (1)$: By Theorems \ref{dog_md}, and \ref{dog for sum}, it is clear.   \\
$ (3) \Longrightarrow (1)$: It is obvious.\\
 
 $(1) \Longrightarrow (3)$: Define  $ S, T \in L ( E \oplus E) $
\[
 S =  \begin{bmatrix}
 A  &  0  \\ 
0 & 0
\end{bmatrix}, \qquad 
 T =  \begin{bmatrix}
 B_{1} &  B_{2}  \\ 
0 & 0
\end{bmatrix}  \].
 Then $ R(S) \subset R(T) $, by Theorem \ref{dog_md}, suppose 
  \[
  X = 
 \begin{bmatrix}
 X_{1}  &  X_{3} \\ 
	 X_{2}  &  X_{4}
 \end{bmatrix}   
\]
is the solution of $ S = TX $, so we have $ A = B_{1}X_{1} + B_{2}X_{2}. $ This completes the proof.
\end{proof}
\end{theorem}

Now we want to show that under certain conditions the sum of two $ { K-\text{frame}} $, is a $K$-frame. Firstly   	
	suppose $ \{x_{n}\} _{{n \in \mathbb{J} }}$ and $ \{y_{n}\} _{{n \in \mathbb{J}}} $ are two Bessel sequences in Hilbert module $E$, then by the Minkowski's inequality it is easy to see that $ \{x_{n} + y_{n} \} _{{n \in \mathbb{J}}} $ is also a Bessel sequence for $E.$

 \begin{theorem}
 Let $ \{x_{n}\} _{{n \in \mathbb{J}}} $ and $ \{y_{n}\} _{{n \in \mathbb{J}}} $ be two $K$-frames for $E$ and also let the corresponding operators in Theorem \ref{L oper}, be $ L_{1} $ and $ L_{2} $ respectively. If $ L_{1}{L_{2}}^{*} $ and $ L_{2}{L_{1}}^{*} $ are positive operators and
 $ R(L_{1}) + R(L_{2}) $ is closed, then $ \{ x_{n} + y_{n} \} _{{n \in \mathbb{J}}} $ is a $K$-frame for $E$.
 \begin{proof}
By the hypothesis we have $$ L_{1}e_{n} = x_{n},  L_{2}e_{n} = y_{n}, R(K) \subset R(L_{1}), R(K) \subset R(L_{2}) $$, where 
  $ \{e_{n}\} _{{n \in \mathbb{J}}} $ is the canonical orthonormal basis of ${\ell}^{2}(\mathcal{A})$. So  $ R(K) \subset R(L_{1}) + R(L_{2}) $, by Theorem \ref{dog for sum}, $ KK^{*} \leq \lambda ( L_{1}{L_{1}}^{*} + L_{2}{L_{2}}^{*} ) $ for some $ \lambda > 0 $. On the other hand for each $x \in E $, 
  \begin{equation*}
	\begin{split}
	 \sum_{n=1}^\infty  \langle x , x_{n} + y_{n} \rangle\langle  x_{n} + y_{n} , x \rangle & = \sum_{n=1}^\infty  \langle ( {L_{1}}^{*} + {L_{2}}^{*}) x , e_{n}  \rangle\langle  e_{n} , {L_{1}}^{*} + {L_{2}}^{*}) x \rangle \\
	 & = \sum_{n=1}^\infty  \langle ({L_{1}} + {L_{2}})^{*} x , e_{n}  \rangle\langle  e_{n} , ({L_{1}} + {L_{2}})^{*} x \rangle \\
	 & = \Vert ({L_{1}} + {L_{2}})^{*} x \Vert_{l^{2}(\mathcal{A})} \\
	 & = \langle ({L_{1}} + {L_{2}})^{*} x , ({L_{1}} + {L_{2}})^{*} x \rangle \\
	 & = \langle {L_{1}}^{*} x , {L_{1}}^{*} x \rangle + \langle {L_{1}}^{*} x , {L_{2}}^{*} x \rangle \\
	 & + \langle {L_{2}}^{*} x , {L_{1}}^{*} x \rangle +
	 \langle {L_{2}}^{*} x , {L_{2}}^{*} x \rangle \\
	 & \geqslant \langle (L_{1}{L_{1}}^{*} + L_{2}{L_{2}}^{*}) x ,  x \rangle \\
	 & \geqslant \frac{1}{\lambda} ( \langle KK^{*}x , x \rangle \\
	 & \geqslant \frac{1}{\lambda} ( \langle K^{*}x , K^{*}x \rangle .
	 \end{split}
 \end{equation*}
Thus $ \{ x_{n} + y_{n} \} _{{n \in \mathbb{J}}} $ is a $ K- $frame.
\end{proof}  
 
 \end{theorem}
\section{$K$-frame vectors for unitary systems}
	A unitary system is a set of unitary operators contains the identity operator.
	A vector $ \psi $ in $E$ is called a \emph{ complete K-frame } vector for a unitary system $\mathcal{U}$ on $E$ if $\mathcal{U}\psi = \{U\psi \mid U \in \mathcal{U} \} $ is a K-frame for $E$. 
	If $\mathcal{U}\psi$ is an orthonormal basis for E, then $ \psi $ is called a \emph{ complete wandering} vector for $\mathcal{U}$.
	The set of all complete K-frame vectors and complete wandering vectors for $\mathcal{U}$ is denoted by $\mathcal{F}_{K}(\mathcal{U}) $ and $ \omega(\mathcal{U})$, respectively.
	In this section we characterize $\mathcal{F}_{K}(\mathcal{U}) $ in terms of operators and elements of $ \omega(\mathcal{U}) $. Also we give conditions under which a linear operation on given elements of $\mathcal{F}_{K}(\mathcal{U}) $ remain an element of $\mathcal{F}_{K}(\mathcal{U}) $.
	\begin{definition}
		For unitary system $ \mathcal{U} $ on Hilbert module $E$ and $ \psi \in E $, the local commutant $\mathcal{C}_{\psi}(\mathcal{U}) $ of $ \mathcal{U} $ at $ \psi $ is defined by 
		$$\mathcal{C}_{\psi}(\mathcal{U}) = \{ T \in L(E) \mid TU\psi = UT\psi, \quad  U \in \mathcal{U}  \}.$$
		Also let $ {\ell^{2}_{\mathcal{U}}}(\mathcal{A}) $  be the Hilbert $\mathcal{A}$-module defined by
		$$ {\ell^{2}_\mathcal{U}}(\mathcal{A}) = \{  \{ a_{U} \} \subset \mathcal{A} \quad : \sum a_{U} a_{U}^{*} \quad \text{converges in} \quad \lVert . \lVert \}. $$		
	\end{definition}

The following theorem characterizes complete $ K $-frame vectors in terms of operators on  complete wandering vectors.
\begin{theorem}
	Suppose $  \mathcal{U} $ is a unitary system of $ E $, $ K \in L(E) $, \linebreak $ \psi \in \omega(\mathcal{U})   $, $ \eta \in E $, and suppose that $ \psi_{\eta} \in L(E , {\ell^{2}_{\mathcal{U}}}(\mathcal{A}) )$ is given by $ T_{\eta}(x) = \{ \langle x , U_{\eta} \rangle \}_{ U \in {\mathcal{U}} }$ and $ \overline{ R({T_{\eta}}^{*}) } $ is orthogonal complemented.
	   Then $  \eta \in \mathcal{F}_{K}(\mathcal{U})  $ if and only if there exist an operator $ A \in \mathcal{C}_{\psi}(\mathcal{U}) $ with $ R(K) \subset R(A) $ such  that $ \eta = A\psi $.
	
	\begin{proof}
	
	($ \Longrightarrow $) Suppose $ \{ e_{U} \}_{U \in \mathcal{U}} $ denote the standard orthonormal basis of $ {\ell^{2}_{\mathcal{U}}}(\mathcal{A}) $, where $ e_{U} $ takes value $ 1_{ \mathcal{A} }$ at $ U $ and $ 0_{\mathcal{A}} $ at every where  else.
	Now suppose   $  \eta \in \mathcal{F}_{K}(\mathcal{U})  $, define operator $ T_{\psi} $ from $ E $ to $ {\ell^{2}_{\mathcal{U}}}(\mathcal{A}) $ by
	  $ {T_{\psi} x = \sum_{U \in \mathcal{U} } \langle x,U_{\psi}\rangle e_{U}}  $.	It is easy to check that $ T_{\psi}$ is well defined, adjointable and invertible. Let $  A={T_{\eta}}^{*} T_{\psi} $. Then for any $ x\in E $, we have 
	$ 	Ax =\sum_{U \in \mathcal{U}} \langle x,U_{\psi}\rangle U_{\eta} $ and  $  A^{*}x = \sum_{U \in \mathcal{U}} \langle x,U_{\eta}\rangle U_{\psi} $, also 
		
		\begin{equation}
	\begin{split}
		\langle A^{*}x , A^{*}x \rangle & =\langle \sum_{U \in \mathcal{U} } \langle x , U\eta \rangle U\psi , \sum_{U \in \mathcal{U} } \langle x , U\eta \rangle U\psi \rangle \\
	& = \sum_{U \in \mathcal{U} } \langle x, U\eta \rangle \langle U\eta , x \rangle \\
	& \geq c \langle K^{*}x , K^{*}x \rangle,
	\end{split}
	\end{equation}
	
	where $ c > 0 $ is the lower bound for $K$-frame 	$ \{ U\eta \mid U \in \mathcal{U} \} $.
On the other hand $ R(A) = R({T_{\eta}}^{*})$ and so by Theorem \ref*{dog_md}, we have $ R(K) \subset R(A) $. To complete the proof, it remains to prove that $ \eta = A\psi $ and $ A \in \mathcal{C}_{\psi}(\mathcal{U}) $.
For any $ U $ and $ V $ in $ \mathcal{U} $\\
  \begin{equation}
  \begin{split}
  \langle V\eta , AU\psi \rangle & = \langle V\eta , \sum_{U \in \mathcal{U} } \langle U\psi , W\psi \rangle W\eta \rangle \\
  & = \sum_{U \in \mathcal{U} } \langle V\eta ,  W\eta \rangle \langle W\psi , U\psi  \rangle \\
  & = \langle V\psi , U\psi \rangle.
  \end{split}
  \end{equation}
This implies that $ AU\psi =U\eta $, so $ A\psi =\eta $.
Also $ AU\psi =U\eta = UA\psi $, hence $  A \in \mathcal{C}_{\psi}(\mathcal{U}) $   and this completes the proof of this part.\\

($ \Longleftarrow $): Suppose that there exists an operator
$  A \in \mathcal{C}_{\psi}(\mathcal{U}) $  with $ R(K) \subset R(A) $ such that $ \eta = A\psi $. Then  for any $ x \in E $, we have\\
 \begin{equation}
 \begin{split}
  \sum_{U \in \mathcal{U} } \langle x , U\eta \rangle \langle U\eta , x \rangle & = \sum_{U \in \mathcal{U} } \langle x , UA\psi \rangle \langle UA\psi , x \rangle \\
  & = \sum_{U \in \mathcal{U} } \langle A^{*}x , U\psi \rangle \langle U\psi , A^{*}x \rangle \\
  & = \langle A^{*}x , A^{*}x  \rangle \\
  & \leq \lVert A^{*} \rVert^{2} \lVert x \rVert ^{2}
  \end{split}
 \end{equation}
So $ \{ U\eta \mid U\in \mathcal{U} \} $ is a Bessel sequence for E.
Now let $ T_{\eta} $  and $ T_{\psi} $ be the operators as we defined in the first part of the proof, since $ \eta = A\psi $ so we have $ T_{\eta} = T_{\psi}A^{*} $. Since $ \psi \in w(\mathcal{U}) $, it is easy to see that $ T^{*}_{\psi} $  is invertible  and hence $ R(T^{*}_{\eta} )=R(A) $. So $ R(K) \subset R(T^{*}_{\eta} ) $. Therefore by using Theorem 3.2 of [10] $ \eta \in \mathcal{U}_{K}(\mathcal{U}) $.

	\end{proof}
\end{theorem}

\bibliographystyle{amsplain}


\end{document}